\newcommand{\switch}[2]{#1} 
\documentclass{article}
\usepackage{quad}

\hypersetup{pdftitle={Five-point Toponogov theorem},
pdfauthor={Nina Lebedeva and Anton Petrunin}}

\newcommand{\Addresses}{{\bigskip\footnotesize

\noindent Nina Lebedeva,
\par\nopagebreak
 \textsc{Saint Petersburg State University, 7/9 Universitetskaya nab., St. Petersburg, 199034, Russia}
\par
\nopagebreak
 \textsc{St. Petersburg Department of V.A. Steklov Institute of Mathematics of the Russian Academy of Sciences, 27 Fontanka nab., St. Petersburg, 191023, Russia}
  \par\nopagebreak
  \textit{Email}: \texttt{lebed@pdmi.ras.ru}

\medskip

\noindent   Anton Petrunin, 
\par\nopagebreak
 \textsc{Math. Dept. PSU, University Park, PA 16802, USA.}
  \par\nopagebreak
  \textit{Email}: \texttt{petrunin@math.psu.edu}
  
}}

\begin{document}

\title{Five-point Toponogov theorem}
\author{Nina Lebedeva and Anton Petrunin}

\date{}
\maketitle
\begin{abstract}
We give an if-and-only-if condition on five-point metric spaces that admit isometric embeddings into complete nonnegatively curved Riemannian manifolds.
\end{abstract}

\section{Introduction}

Toponogov theorem provides an if-and-only-if condition on a metric on four-point space that admits an isometric embedding into a complete nonnegatively curved Riemannian manifold.
The only-if part is proved by Victor Toponogov, and the if part follows from a result of Abraham Wald \cite[\S 7]{wald}.

We show that the so-called Lang--Schroeder--Sturm inequality is the analogous condition for five-point spaces.
 
The only-if part is well-known, but the if part is new.
It was hard to imagine some new restrictions on five-point sets, but now we know there are none.

Let us formulate the Lang--Schroeder--Sturm inequality.
Consider an $(n+1)$-point array $(p,x_1,\dots x_n)$ in a metric space $X$.
We say that the array satisfies Lang--Schroeder--Sturm inequality with center $p$
if for any nonnegative values $\lambda_1,\dots,\lambda_n$ we have
\[\sum_{i,j}a_{ij}\cdot \lambda_i\cdot\lambda_j\ge 0,\]
where $a_{ij}=|p-x_i|_X^2+|p-x_j|_X^2-|x_i-x_j|_X^2$
and we denote by $|\ -\ |_X$ the distance between points in $X$.

Recall that any point array in a  complete nonnegatively curved Riemannian manifold (and, more generally, in any nonnegatively curved Alexandrov space) meets the Lang--Schroeder--Sturm inequality \cite{lang-schroeder, sturm}.
In particular, 
the Lang--Schroeder--Sturm inequalities for all relabelings of points in a finite metric space $F$
gives a necessary condition for the existence of isometric embedding of $F$ into a complete Riemannian manifold with nonnegative curvature.
In this note, we show that this condition is sufficient if $F$ has at most 5 points.

\begin{thm}{Theorem}\label{thm:main}
A five-point metric space $F$ admits an isometric embedding into a complete nonnegatively curved Riemannian manifold
if and only if all Lang--Schroeder--Sturm inequalities hold in $F$.
\end{thm}

\begin{wrapfigure}{r}{15mm}
\vskip-0mm
\centering
\includegraphics{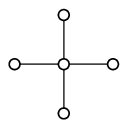}
\end{wrapfigure}

In the next section, we will give a reformulation of the theorem using the so-called \emph{(4+1)-point comparison} \cite{alexander2019alexandrov,AKP-Kirszbraun} which is also equivalent to \emph{graph comparison} \cite{lebedeva-petrunin-zolotov} for the star graph shown on the diagram.

Since we know that Lang--Schroeder--Sturm inequalities are necessary,
it remains to construct a complete nonnegatively curved Riemannian manifold that contains an isometric copy of a given 5-point space satisfying the assumptions.

Our proof uses a brute-force search of certain configurations that was originally done on a computer.
We present a hand-made proof that was found later.
It is still based on brute-force search, and we hope that a more conceptual proof will be found.
Our paper is inspired by the note of Vladimir Zolotov and the first author \cite{lebedeva-zolotov};
the results from this note are discussed briefly in the last section.

\parbf{Acknowledgments.}
We want to thank Arseniy Akopyan and Alexander Gil for helping us with programming.
We would also like to thank Tadashi Fujioka and Tetsu Toyoda for pointing out errors and misprints in the preliminary version of this paper.
The first author was partially supported by Russian Foundation for Basic Research grant 20-01-00070; 
the second author was partially supported by National Science Foundation grant DMS-2005279.

\section{LSS(\textit{n}) and (\textit{n}+1)-comparison}

The (\textit{n}+1)-comparison is another condition that holds for any $(n+1)$-point array in  Alexandrov spaces \cite{alexander2019alexandrov,AKP-Kirszbraun}.
It says that given a point array $p,x_1,\dots,x_n$ in a nonnegatively curved Alexandrov space $A$ 
there is an array $\tilde p,\tilde x_1,\dots,\tilde x_n$ in a Hilbert space $\HH$ such that 
\[
|\tilde p-\tilde x_i|_{\HH}=|p-x_i|_{A}
\quad
\text{and}
\quad
|\tilde x_i-\tilde x_j|_{\HH}\ge |x_i-x_j|_{A}.
\]
for all $i$ and $j$.
Point $p$ will be called the \emph{center} of comparison.

Let us denote by $S_n$ the star graph of order $n$;
one central vertex in $S_n$ is connected to the remaining $n$.
It is easy to see that (\textit{n}+1)-comparison is equivalent to the $S_n$-comparison --- a particular type of graph comparison introduced in \cite{lebedeva-petrunin-zolotov}.

For general metric spaces, the (\textit{n}+1)-comparison implies the Lang--Schroeder--Sturm inequality,  briefly $\LSS(n)$.
For $n\ge 5$ the converse does not hold \cite[Section 8]{lebedeva-petrunin-zolotov}.
In this section, we will show that these two conditions are equivalent for $n\le 4$.

\begin{thm}{Claim}\label{clm:(4+1)=LSS(4)}
For any 5-point array $p,x_1,\dots,x_4$, the
$\LSS(4)$-inequality is equivalent to the $(4+1)$-comparison.
\end{thm}

Applying the claim, we get the following reformulation of the main theorem.

\begin{thm}{Reformulation}\label{thm:main-(4+1)}
A five-point metric space $F$ admits an isometric embedding into a complete Riemannian manifold with nonnegative curvature
if and only if it satisfies (4+1)-comparison for all relabelings.
\end{thm}

The following proof is nearly identical to the proof of Proposition 4.1 in \cite{lebedeva-petrunin-zolotov}.

\parit{Proof of \ref{clm:(4+1)=LSS(4)}.}
Suppose $p,x_1,\dots,x_4$ satisfies $\LSS(4)$;
we need to show that it also meets the (4+1)-comparison.

Choose a smooth function $\phi\:\RR\to \RR$ such that $\phi(x)=0$ if $x\ge0$ and $\phi(x)>0$, $\phi'(x)<0$ if $x<0$.
Consider a configuration of points $\tilde p,\tilde x_1,\dots,\tilde x_4\in \HH$ such that 
\[
|\tilde p-\tilde x_i|_{\HH}=|p-x_i|_{A}
\]
and the  value
\[
s=\sum_{i<j}\phi(|\tilde x_i-\tilde x_j|_{\HH}- |x_i-x_j|_{A})
\]
is minimal.
Note that $s\ge0$;
if $s=0$, then we get the required configuration.

Suppose $s>0$.
Consider the graph $\Gamma$ with 4 vertices labeled by $\tilde x_1,\tilde x_2,\tilde x_3,\tilde x_4$ such that 
$(\tilde x_i,\tilde x_j)$ is an edge if and only if $|\tilde x_i-\tilde x_j|_{\HH}<|x_i-x_j|_{A}$.
Assume $\tilde x_1$ has a single incident edge, say $(\tilde x_1,\tilde x_2)$.
Since $s$ takes minimal value, we have 
\begin{align*}
|\tilde x_1-\tilde x_2|_{\HH}&=|\tilde x_1-\tilde p|_{\HH}+|\tilde x_2-\tilde p|_{\HH}=
\\
&=| x_1- p|_{A}+| x_2- p|_{A}\ge
\\
&\ge |x_i-x_j|_{A}
\end{align*}
--- a contradiction.
It follows that $\Gamma$ contains no end-vertices.
Therefore it is isomorphic to one of the four graphs on the diagram.

\begin{figure}[ht!]
\centering
\includegraphics{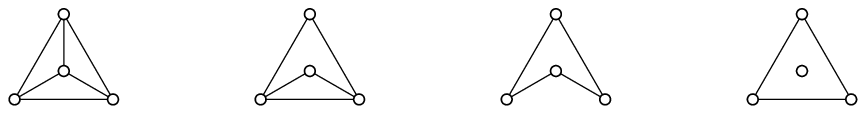}
\end{figure}

Without loss of generality, we can assume that $\tilde p=0$.
Note that any point $\tilde x_i$ cannot lie in an open half-space with all its
adjacent points.
Indeed, assume it does, then $\tilde x_i$ and all its adjacent points lie in a finite-dimensional half-space;
denote by $\Pi$ its boundary plane.
Then rotating $\tilde x_i$  slightly around $\Pi$ will increase the distances from $\tilde x_i$ to all its adjacent vertices;
so the value $s$ will decrease --- a contradiction.

In the 6-edge case, $\tilde p=0$ lies in the convex hull of $\{\tilde x_1,\tilde x_2,\tilde x_3,\tilde x_4\}$.
In particular, $0=\lambda_1\cdot \tilde x_1+\dots+\lambda_4\cdot\tilde x_4$ for some $\lambda_i\ge0$ such that $\lambda_1+\lambda_2+\lambda_3+\lambda_4=1$.
The latter contradicts $\LSS(4)$.

Similarly, in the 5- and 3-edge cases, we can assume that $\tilde x_1\tilde x_2\tilde x_3$ is a 3-cycle of $\Gamma$.
In this case, $0=\lambda_1\cdot \tilde x_1+\lambda_2\cdot\tilde x_2+\lambda_3\cdot\tilde x_3$ for some $\lambda_i\ge0$ such that $\lambda_1+\lambda_2+\lambda_3=1$, and we arrive at a contradiction with $\LSS(3)$.

Finally, the 4-edge graph (the 4-cycle) cannot occur.
In this case, we may think that $\tilde x_1,\tilde x_2,\tilde x_3,\tilde x_4$ is the 4-cycle.
Note that the points $\tilde x_1,\tilde x_2,\tilde x_3,\tilde x_4$  lie in one plane so that the direction of $\tilde x_1$ is opposite to $\tilde x_3$,
and the direction of $\tilde x_2$ is opposite to $\tilde x_4$.
Let us think that this is the horizontal plane in $\RR^3$.
Then rotating the pair $\tilde x_1$, $\tilde x_3$ slightly up and 
the pair $\tilde x_2$, $\tilde x_4$  slightly down, decreases~$s$ --- a contradiction.
\qeds

\section{Associated form}

In this section, we recall a construction from \cite{petrunin-2017}.
Let $\bm{x}=(x_1,\dots,x_n)$ be a point array in a metric space $X$.

Choose a simplex $\triangle$ in $\RR^{n-1}$; for example, we can take the standard simplex with the first $(n-1)$ of its vertices $v_1,\dots,v_{n-1}$ form the standard basis on $\RR^{n-1}$, and $v_n=0$.

Consider a quadratic form $W_{\bm{x}}$ on $\RR^{n-1}$ that is uniquely defined by
\[W_{\bm{x}}(v_i-v_j)=|x_i-x_j|^2_X\] 
for all $i$ and $j$.
It will be called
the \emph{associated} form to the point array $\bm{x}$.
The following claim is self-evident:

\begin{thm}{Claim}\label{clm:W>=0}
An array $\bm{x}=(x_1,\dots,x_n)$ with a (semi)metric is isometric to an array in a Euclidean space if and only if 
$W_{\bm{x}}(v)\ge 0$
for any $v\in \RR^{n-1}$.
\end{thm}

In particular, the condition $W_{\bm{x}}\ge 0$ for a triple $\bm{x}=(x_1,x_2,x_3)$ is equivalent to 
the three triangle inequalities for the distances between $x_1$, $x_2$, and $x_3$.
For an $n$-point array, it implies that $W_{\bm{x}}(v)\ge 0$ for any vector $v$ in a plane spanned by a triple of vertices of $\triangle$.

\parbf{Lang--Schroeder--Sturm inequalities.}
Consider lines that connect a point on a facet of $\triangle$ with its opposite vertex.
The union of these lines forms a cone in $\RR^{n-1}$; denote it by $K_n$.
Note that $K_3=\RR^2$, but for $n\ge 4$ the cone $K_n$ is a proper subset of $\RR^{n-1}$.

The following claim is a reformulation of Lang--Schroeder--Sturm inequalities for all relabeling of $\bm{x}$:

\begin{thm}{Claim}\label{clm:W(w)<0}
Let $\bm{x}=(x_1,\dots,x_n)$ be an $n$-point array in an Alexandrov space with nonnegative curvature.
Then $W_{\bm{x}}(w)\ge0$ for any $w\in K_n$.
\end{thm}

\section{Tense arrays}

Assume $(p,x_1,\dots,x_n)$ is an array of points with a metric that satisfies $\LSS(n)$ with center at $p$.
Suppose $n\le 4$; by Claim~\ref{clm:(4+1)=LSS(4)}, we have a comparison configuration 
$(\tilde p,\tilde x_1,\dots,\tilde x_n)$ with center $p$.

We say that an array $(p,x_1,\dots,x_n)$ is \emph{tense} with center $p$ if the comparison configuration $(\tilde p,\tilde x_1,\dots,\tilde x_n)$ is unique up to congruence and isometric to the original array.

Note that if $(p,x_1,\dots,x_n)$ is tense, then in its comparison configuration $\tilde p$ lies in the convex hull of the remaining points $\tilde x_1,\dots,\tilde x_n$.
In particular, 
\[\sum_i\lambda_i\cdot (\tilde x_i-\tilde p)=0 \leqno({*})\]
for some nonnegative coefficients $\lambda_1,\dots, \lambda_n$ such that $\lambda_1+\dots+\lambda_n=1$.
If we can choose all positive $\lambda_i$ in $({*})$, then we say that $(p,x_1,\dots,x_n)$ is a \emph{nondegenerate tense array}.
The following statement describes nondegenerate tense arrays.

\begin{thm}{Claim}\label{clm:nondeg-tense}
Assume $n\le 4$ and $(p,x_1,\dots,x_n)$ is an array of points with a metric that satisfies $\LSS(n)$ with center $p$.
Suppose we have equality in $\LSS(n)$ for some positive $\lambda$-parameters;
that is, for some positive values $\lambda_1,\dots,\lambda_n$ we have
\[\sum_{i,j}a_{ij}\cdot \lambda_i\cdot\lambda_j= 0,\]
where $a_{ij}=|p-x_i|^2+|p-x_j|^2-|p-x_j|^2$.
Then $(p,x_1,\dots,x_n)$ is a nondegenerate tense array with center $p$.
\end{thm}

\parit{Proof.}
By Claim~\ref{clm:(4+1)=LSS(4)}, we have a comparison configuration $(\tilde p,\tilde x_1,\dots,\tilde x_n)$ with center $p$.
Set 
\[\tilde a_{ij}
=
\langle\tilde x_i-\tilde p,\tilde x_j-\tilde p\rangle
=
|\tilde p-\tilde x_i|^2+|\tilde p-\tilde x_j|^2-|\tilde x_i-\tilde x_j|^2.\]
Since $|\tilde x_i-\tilde x_j|\ge| x_i- x_j|$, we have $\tilde a_{ij}\le a_{ij}$ for any $i$ and $j$.
Note that 
\begin{align*}
0&=\sum_{i,j}a_{ij}\cdot \lambda_i\cdot\lambda_j\ge
\sum_{i,j}\tilde a_{ij}\cdot \lambda_i\cdot\lambda_j= 
2\cdot\left|\sum_i\lambda_i\cdot (\tilde x_i-\tilde p)\right|^2\ge 
0.
\end{align*}
It follows that 
\[\sum_i\lambda_i\cdot (\tilde x_i-\tilde p)=0.\]
Further, $\tilde a_{ij}= a_{ij}$, and, therefore, $|\tilde x_i-\tilde x_j|=| x_i- x_j|$ for all $i$ and $j$ --- hence the result.
\qeds

Note that any 2-point array is a degenerate tense array; the center can be chosen arbitrarily.

A 3-point array $(p,x,y)$ is tense with center $p$ if we have equality 
\[|p-x|+|p-y|=|x-y|.\]
Note that any tense 3-point array with distinct points is nondegenerate.

Let $(p,x,y,z)$ be a tense 4-point array.
Then there is an isometric comparison configuration $(\tilde p,\tilde x,\tilde y,\tilde z)$ with $\tilde p$ lying in the solid triangle $\tilde x\tilde y\tilde z$.
Suppose all points $p$, $x$, $y$, and $z$ are distinct.
If $\tilde p$ lies in the interior of the solid triangle or the triangle is degenerate, then the array $(\tilde p,\tilde x,\tilde y,\tilde z)$ is nondegenerate.
Otherwise, if $\tilde p$ lies on a side, say $[\tilde x,\tilde y]$, and the triangle $[\tilde x\tilde y\tilde z]$ is nondegenerate, then $(\tilde p,\tilde x,\tilde y,\tilde z)$ is degenerate.
In the latter case, the 3-point array $(\tilde p,\tilde x,\tilde y)$ is tense and nondegenerate.

\begin{thm}{Claim}\label{clm:10-2k}
Let $\bm{x}=(x_1,\dots,x_5)$ be a 5-point array in a metric space that satisfies all $\LSS(4)$-inequalities.
Suppose that $\bm{x}$ has $k$ three-point tense arrays and no tense arrays with four and five points.
If $k\le 4$ then there is a $2$-dimensional subspace $S$ of quadratic forms on $\RR^4$ such that for any form $U\in S$ that is sufficiently close to zero the array with associated form $W_{\bm{x}}+U$ satisfies all $\LSS(4)$-inequalities.
\end{thm}

\parit{Proof.}
Given a tense three-point subarray of $\bm{x}$, say $(x_1,x_2,x_3)$, consider the metrics on $\bm{x}$ 
such that the three distances between $x_1$, $x_2$, and $x_3$ are proportional to the original distances, and the remaining distances are arbitrary.
This set defines a subspace of quadratic forms of codimension 2 in the 10-dimensional space of quadratic forms on $\RR^4$.
Since $k\le 4$, 
taking the intersection of all such subspaces we get a subspace $S$ of dimension at least $2$.

It remains to show that $S$ meets the claim --- assume not.
That is, for arbitrary small $U\in S$ the metric on $\bm{x}$ defined by the associated quadratic form $W_{\bm{x}}+U$ does not satisfy the $\LSS(4)$-inequality.
It means that there is a vector $w\in K_5$ such that 
\[W_{\bm{x}}(w)+U(w)< 0.\]

Choose a positive quadratic form $I$ and minimal $t>0$ such that  
\[W_{\bm{x}}(w)+U(w)+t\cdot I(w)\ge 0
\leqno({*}{*})\] 
for any $w\in K_5$.

Note that for some $w\in K_5$, we have equality in $({*}{*})$.
The metric on $\bm{x}$ that corresponds to the form $W_{\bm{x}}+U+t\cdot I$ satisfies all $\LSS(4)$ inequalities and
by \ref{clm:nondeg-tense}, it has a tense array with at least 3 points.

Choose an array $Q$ that remains to be tense as $U\to 0$. 
Note that $Q$ is isometric to an array in Euclidean space.
Since $W_{\bm{x}}+U+t\cdot I\to W_{\bm{x}}$ as $U\to 0$, the array $Q$ must contain one of the three-point tense arrays for the original metric.
The latter is impossible since $t>0$
--- a contradiction.
\qeds

\section{Extremal metrics}\label{sec:ext}

Denote by $\mathcal{A}_5$ the space of metrics on a 5-point set $F=\{a,b,c,d,e\}$ that admits an embedding into a Riemannian manifold with nonnegative curvature.
The associated quadratic forms for spaces in $\mathcal{A}_5$ form a convex cone in the space of all quadratic forms on $\RR^4$.
The latter follows since nonnegative curvature survives after rescaling and passing to a product space.

Denote by $\mathcal{B}_5$ the space of metrics on $F$ that satisfies all Lang--Schroeder--Sturm inequalities for all relabelings.
As well as for $\mathcal{A}_5$, the associated forms for spaces in $\mathcal{B}_5$ form a convex cone in the space of all quadratic forms on $\RR^4$.

Since the associated quadratic form describes its metric completely, we may identify $\mathcal{A}_5$ and $\mathcal{B}_5$ with subsets in $\RR^{10}$ --- the space of quadratic forms on $\RR^4$.
This way we can think that $\mathcal{A}_5$ and $\mathcal{B}_5$ are convex cones in $\RR^{10}$.

The set $\mathcal{B}_5$ is a cone so it does not have extremal points except the origin.
The origin corresponds to degenerate metric with all zero distances.
But $\mathcal{B}_5$ is a cone over a convex compact set $\mathcal{B}_5'$ in the sphere $\mathbb{S}^9\subset \RR^{10}$.
The extremal points of $\mathcal{B}_5'$ correspond to extremal rays of $\mathcal{B}_5$;
metrics on extremal rays will be called \emph{extremal}.
Note that if an extremal metric $\rho$ lies in the interior of a line segment between metrics $\rho'$ and $\rho''$ in $\mathcal{B}_5$, then both metrics $\rho'$ and $\rho''$ are proportional to $\rho$.

Since Lang--Schroeder--Sturm inequalities are necessary for the existence of isometric embedding into a complete nonnegatively curved Riemannian manifold,
we have that 
\[\mathcal{A}_5\subset\mathcal{B}_5.\]
To prove the theorem we need to show that the opposite inclusion holds as well.
Since $\mathcal{B}_5$ is the convex hull of its extremal metrics, it is sufficient to prove the following:

\begin{thm}{Proposition}\label{prop:main}
Given an extremal space $F$ in $\mathcal{B}_5$, there is a complete nonnegatively curved Riemannian manifold that contains an isometric copy of $F$.
\end{thm}

\parit{Proof.}
Note that any extremal space  $F$ contains a tense set.
If not, then an arbitrary slight change of metric keeps it in $\mathcal{B}_5$ which is impossible for an extremal metric.

The remaining part of the proof is broken into cases:
\begin{itemize}
\item $F$ contains a 5-point tense set.
In this case, $F$ admits an isometric embedding into Euclidean space --- the problem is solved.
\item $F$ contains a 4-point  tense set.
This case follows from Proposition~\ref{prop:4-tense} below.
\item $F$ contains only 3-point  tense sets. This is the hardest part of the proof; it follows from Proposition~\ref{prop:3-tense}.
\qeds
\end{itemize}

\section{Four-point tense set}\label{sec:4-tense}

\begin{thm}{Proposition}\label{prop:4-tense}
Suppose that a 5-point metric space $F$ satisfies all Lang--Schroeder--Sturm inequalities and contains a 4-point tense set.
Then $F$ is isometric to a subset of а complete nonnegatively curved Riemannian manifold. 
\end{thm}

In the following proof, we first construct a nonnegatively curved Alexandrov space with an isometric copy of $F$ and then smooth it.
The space will be a doubling of a convex polyhedral set in $\RR^3$.

We use notations
\[\hinge xyz,\quad\measuredangle\hinge xyz,\quad\text{and}\quad \angk{x}{y}{z}\]
for hinge, its angle measure, and the model angle respectively.

\parit{Proof.}
Let us label the points in $F$ by $p$, $q$, $x_1$, $x_2$, and $x_3$ so that the array $(p,x_1,x_2,x_3)$ is tense with center $p$.

By the definition of a tense array, we can choose an array $(\tilde p, \tilde x_1, \tilde x_2, \tilde x_3)$ in $\RR^2$ that is isometric to  $(p, x_1, x_2, x_3)$.
Consider $\RR^2$ as a plane in $\RR^3$.

By \ref{clm:(4+1)=LSS(4)}, we can apply the (4+1)-comparison.
It implies the existence of point $\tilde q\in\RR^3$ such that
\[
|\tilde p-\tilde q|_{\RR^3}=|p-q|_F
\quad\text{and}\quad
|\tilde x_i-\tilde q|_{\RR^3}\ge|x_i-q|_F
\leqno({*})
\]
for any $i$. 

Further, let us show that there are points $\tilde q_1$, $\tilde q_2$, $\tilde q_3$ in the plane thru $\tilde p$, $\tilde x_1$, $\tilde x_2$, and $\tilde x_3$ such that the following four conditions
\[|\tilde p-\tilde q_i|_{\RR^3}\ge |p-q|_F,
\quad
|\tilde x_j-\tilde q_i|_{\RR^3}\ge|x_j-q|_F,
\quad
|\tilde x_i-\tilde q_i|_{\RR^3}=|x_i-q|_F,
\leqno({*}{*})\]
hold for all $i$ and $j$.

By these conditions, $\tilde q_1$ must lie on the circle with the center at $x_1$ and radius $|x_1-q|_F$.
Denote by $\Gamma_1$ the intersection of this circle with the angle vertical to the hinge $\hinge{\tilde x_1}{\tilde x_2}{\tilde x_3}$.
Let us show that $\tilde q_1$ can be chosen on $\Gamma_1$.
The point $\tilde q_1$ has to satisfy additional three conditions:
\[\measuredangle\hinge{\tilde x_1}{\tilde x_i}{\tilde q_1}\ge \angk{x_1}{x_i}{q},
\qquad
\measuredangle\hinge{\tilde x_1}{\tilde p}{\tilde q_1}\ge \angk{x_1}{p}{q}\]
for $i\ne 1$.
Each condition describes a subarc of $\Gamma_1$, say $\breve{X}_{1,2}$, $\breve{X}_{1,3}$, and $\breve{P}_{1}$.

\begin{wrapfigure}[5]{r}{30mm}
\vskip-2mm
\centering
\includegraphics{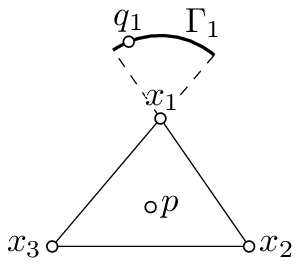}
\end{wrapfigure}

By the construction and comparison we have 
\begin{align*}
\angk{x_1}{q}{x_2}+\angk{x_1}{q}{x_3}+\angk{x_1}{x_2}{x_3}&\le 2\cdot\pi,
\\
\angk{x_1}{q}{p}+\angk{x_1}{q}{x_i}+\angk{x_1}{p}{x_i}&\le 2\cdot\pi,
\\
\angk{x_1}{x_2}{x_3}=\measuredangle\hinge{\tilde x_1}{\tilde x_2}{\tilde x_3},
\qquad
\angk{x_1}{x_i}{p}&=\measuredangle\hinge{\tilde x_1}{\tilde x_i}{\tilde p}
\end{align*}
for $i\ne 1$.
These inequalities and identities imply that each pair of arcs $\breve{X}_{1,2}$, $\breve{X}_{1,3}$, and $\breve{P}_{1}$ have a nonempty intersection.
By 1-dimensional Helly's theorem, all three arcs intersect;
so we can choose $\tilde q_1$ in this intersection.
The same way we construct $\tilde q_2$ and $\tilde q_3$.

Now let us show that there is a point $\tilde s\in\RR^2$ such that
\begin{align*}
|\tilde p-\tilde s|_{\RR^3}&\le|p-q|_F,
&
|\tilde x_i-\tilde s|_{\RR^3}&\le|x_i-q|_F.
\end{align*}
for all $i$.
In other words, the following four closed balls have a nonempty intersection: 
$\bar B[\tilde p,|p-q|_F]$ and
$\bar B[\tilde x_i,|x_i-q|_F]$ for all $i$.
Indeed, by the overlap lemma \cite{alexander2019alexandrov}, any 3 of these balls have a nonempty intersection;
it remains to apply Helly's theorem.
Note that we can assume that $\tilde s$ lies in the convex hull of $\tilde x_1$, $\tilde x_2$, and $\tilde x_3$.

The four perpendicular bisectors to 
$[\tilde s, \tilde q]$, 
$[\tilde s, \tilde q_1]$, 
$[\tilde s, \tilde q_2]$, 
$[\tilde s, \tilde q_3]$ cut from $\RR^3$ a closed convex set $V$ that contains $\tilde s$.
(It might be a one-sided infinite triangular prism or, if $\tilde q$ lies in the plane of the triangle, a two-sided infinite quadrangular prism.)
Note that the inequalities $({*})$ and $({*}{*})$ imply that $V$ contains the points $\tilde p$, $\tilde x_1$, $\tilde x_2$, and $\tilde x_3$.

Consider the doubling $W$ of $V$ with respect to its boundary;
it is an Alexandrov space with nonnegative curvature \cite[5.2]{perelman:spaces2}.
Denote by $\iota_1$ and $\iota_2$ the two isometric embeddings $V\to W$.
By construction, the array $\hat p=\iota_1(\tilde p)$, $\hat x_1=\iota_1(\tilde x_1)$, $\hat x_2=\iota_1(\tilde x_2)$, $\hat x_3=\iota_1(\tilde x_3)$,  $\hat s=\iota_2(\tilde s)$ in $W$ is isometric to the array $(p, q, x_1, x_2, x_3)$ in $F$.

\medskip

Finally, we need to show that the obtained space can be smoothed into a Riemannian manifold that still has an isometric copy of $F$.
This part is divided into two steps; first, we show that the construction above can be made so that the points $\tilde s$, $\tilde p$, $\tilde x_1$, $\tilde x_2$, $\tilde x_3$ do not lie on the edges of $V$.
In this case, there is a compliment, say $U$, of a neighborhood of the singular set in $W$ that contains
the 5-point set together with all the geodesics between them.
After that, we construct a smooth Riemannian manifold with nonnegative curvature that contains an isometric copy of $U$.

\begin{wrapfigure}{r}{30mm}
\vskip-0mm
\centering
\includegraphics{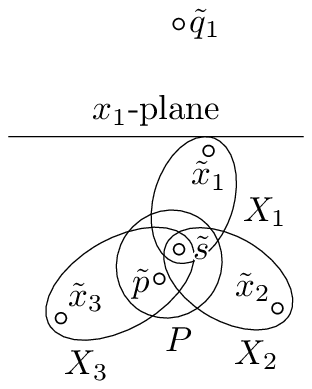}
\end{wrapfigure}

\parit{Step 1.}
Consider the four ellipsoids $P$, $X_1$, $X_2$, $X_3$ with the major axes $|q-p|$, $|q-x_1|$, $|q-x_2|$,$|q-x_3|$, first focus at $\tilde s$, and the second focus at $\tilde p$, $\tilde x_1$, $\tilde x_2$, $\tilde x_3$ respectively.

The construction of the facets of $V$ above implies that each ellipsoid has a tangent plane that contains all the ellipsoids on one side --- these planes are the perpendicular bisectors to $[\tilde s,\tilde q]$ and $[\tilde s,\tilde q_i]$ for all $i$.
Note that any choice of such planes does the trick --- they can be used instead of the perpendicular bisectors discussed above.
These planes will be called \emph{$p$-plane} and \emph{$x_i$-planes} respectively.
We need to choose them so that no pair of these planes pass thru $\tilde p$, $\tilde x_1$, $\tilde x_2$, or $\tilde x_3$.
The latter is only possible if the corresponding ellipsoid degenerates to a line segment.

We can assume that one of the ellipsoids is nondegenerate; otherwise, the array $\tilde p$, $\tilde q$, $\tilde x_1$, $\tilde x_2$, $\tilde x_3$ is isometric to $F$;
in this case, $F$ is isometric to a subset of Euclidean space.
Further, if only one ellipsoid, say $X_1$ is degenerate, then we can move $\tilde s$ slightly making this ellipsoid nondegenerate and keeping the rest of its properties.
So we can assume that three or two ellipsoids are degenerate.

Suppose $X_1$, $X_2$, $X_3$ are degenerate (picture on the left), then it is easy to choose $x_i$-planes tangent to $P$; it solves our problem.
Another triple of ellipsoids, say $P$, $X_1$, $X_2$ might be degenerate only if $\tilde p\in[\tilde x_1,\tilde x_2]$.
This case is even simpler --- we can choose one plane that contains $\tilde p$, $\tilde x_1$, $\tilde x_2$, and tangent to $X_3$.

\begin{figure}[!ht]
\begin{minipage}{.3\textwidth}
\centering
\includegraphics{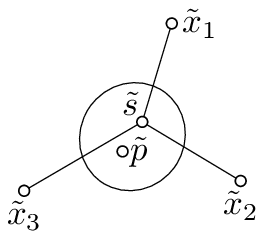}
\end{minipage}
\hfill
\begin{minipage}{.3\textwidth}
\centering
\includegraphics{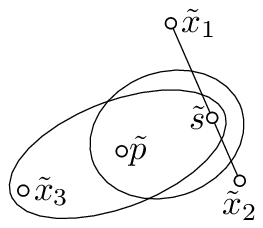}
\end{minipage}
\hfill
\begin{minipage}{.3\textwidth}
\centering
\includegraphics{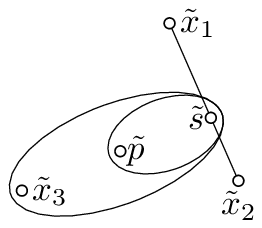}
\end{minipage}
\end{figure}

Now, suppose exactly two ellipsoids are degenerate; note that in this case $P$ is nondegenerate.
Therefore we can assume that $X_1$ and $X_2$ are degenerate.
Further, we can assume that $\tilde s\in[\tilde x_1,\tilde x_2]$;
if not we can slightly move $\tilde s$ toward $\tilde x_1$ and $\tilde x_2$ making $X_1$ and $X_2$ nondegenerate and keep the rest properties of~$\tilde s$.
Since a focus of $P$ lies on $[\tilde x_1,\tilde x_2]$, we have that $x_1$-plane cannot be $x_2$-plane and the other way around.

Suppose that $P$ does not lie in the convex hull of the remaining three ellipsoids
and the same holds for $X_3$ (middle picture).
Then it is easy to make the required choice of planes.

In the remaining case (see picture on the right), either $P$ or $X_3$ lies in the convex hull of the remaining three ellipsoids.
Suppose it is $P$, draw a $p$-plane; note that it is also an $x_3$-plane;
it might be also $x_1$- or $x_2$-plane, but cannot be both.
It remains to add $x_1$-plane and/or $x_2$-plane as needed.
Since an $x_1$-plane cannot be $x_2$-plane and the other way around, we will not get two planes passing thru $x_1$ or $x_2$.

The case when $X_3$ lies in the convex hull of the rest is identical.

\begin{wrapfigure}{r}{40mm}
\vskip-0mm
\centering
\includegraphics{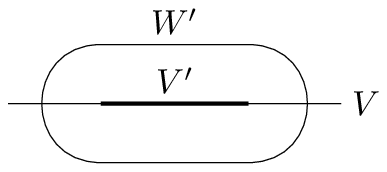}
\end{wrapfigure}

\parit{Step 2.}
Start with the subset $V'\subset V$ that lies on the distance $\pi\cdot\delta$ from its boundary.
Think of $V'$ lying in $\RR^4$, pass to its $2\cdot \delta$-neighborhood.
The boundary of the obtained neighborhood is a convex hypersurface $W'$ in $\RR^4$.
For small $\delta>0$, it meets all the required conditions, except it is only $C^{1,1}$-smooth.
It is straightforward to smooth $W'$ so that the metric changes only near the edges of~$V$.
In this case, the set $F$ remains isometrically embedded in the obtained 3-dimensional manifold.
\qeds

\section{Three-point tense sets}\label{sec:3-tense}

\begin{thm}{Proposition}\label{prop:3-tense}
Suppose that an extremal 5-point metric space $F$ contains only 3-point tense sets.
Then $F$ is isometric to a subset in a nonnegatively curved Riemannian manifold $L$.
Moreover, we can assume that $L$ is homeomorphic to a circle or а plane.
\end{thm}

A three-point tense set $\{a,b,c\}$ with center $b$ will be briefly denoted by $abc$.
Observe that $F$ has tense set $abc$ if and only if 
\[|a-b|_F+|b-c|_F\z=|a-c|_F.\]

{

\begin{wrapfigure}{r}{20mm}
\vskip-6mm
\centering
\includegraphics{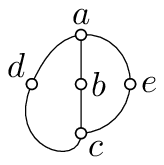}
\end{wrapfigure}

On the diagrams, we will connect three-point tense sets by a smooth curve so that the center is in the middle.
For example, the given diagram corresponds to a metric on $\{a,b,c,d,e\}$ with five tense sets $abc$, $bcd$, $cda$, $dae$, $aec$.

}

\begin{thm}{Classification lemma}\label{lem:key}
Let $F$ be an extremal 5-point metric space; suppose that it has no tense subsets with 4 and 5 points.
Then $F$ has one of three configurations of tense sets shown on the diagram.

In other words, the points in $F$ can be labeled by $\{a,b,c,d,e\}$ so that it has
one of the following three tense-set configurations:

\begin{wrapfigure}{r}{50mm}
\vskip2mm
\centering
\includegraphics{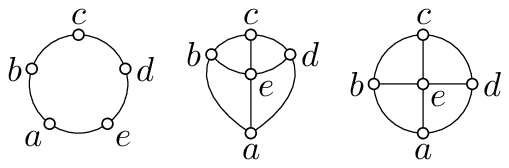}
\end{wrapfigure}
\vskip-6mm
\begin{align*}
&abc, bcd, cde, dea, eab;
\\
&abc, bcd, cda, aec, bed;
\\
&abc, bcd, cda, dab, aec, bed.
\end{align*}

\end{thm}

In the following proof, we use only a small part of this classification.
Namely, we use that it is either the first case (the cycle) or there are two tense sets with a shared center ($bed$ and $cea$).  
However, the proof of this \emph{small part} takes nearly as long as the complete classification.
(We could exclude cases 11, 12 and 16 on page \pageref{pic-210}, but we decided to keep them.) 

\begin{wrapfigure}{r}{20mm}
\vskip-0mm
\centering
\includegraphics{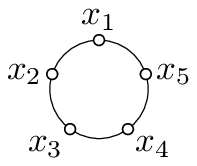}
\end{wrapfigure}

\parit{Proof of \ref{prop:3-tense} modulo \ref{lem:key}.}
Suppose that $F$ has a tense configuration as on the diagram.
In other words, we can label points in $F$  by $\{x_1,x_2,x_3,x_4,x_5\}$ so that
\[|x_{i}-x_{i-1}|_F+|x_{i+1}-x_{i}|_F=|x_{i+1}-x_{i-1}|_F\]
for any $i\pmod 5$.
In this case, $F$ is isometric to a 5-point subset in the circle of length 
$\ell=|x_1-x_2|_F+\dots+|x_4-x_5|_F+|x_5-x_1|_F$.

Now, by the classification lemma, we can assume that two tense triples in $F$ have a common center.
Let us relabel $F$ by $x,v_1,v_2,w_1,w_2$ so that $F$ has tense triples $v_1xv_2$ and $w_1xw_2$.

\begin{wrapfigure}{r}{25mm}
\vskip-6mm
\centering
\includegraphics{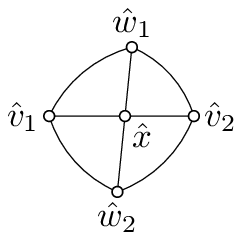}
\end{wrapfigure}

First, we will construct an Alexandrov space $L$ --- a flat disc with at most four singular points.
The disc $L$ will be triangulated by four triangles with vertices $\hat x$, $\hat v_1$, $\hat v_2$, $\hat w_1$, $\hat w_2$ as shown on the diagram.
Each of the four triangles has at most one singular point;
in other words, each triangle is a solid geodesic triangle in a cone.
The sides of the triangles are the same as in $F$.

Note that the metric on the obtained disc is completely determined by the 12 angles of the triangles.
It remains to choose these angles in such a way that $L$ has nonnegative curvature and the map $\iota\:F\to L$ defined by $x\mapsto \hat x$, $v_i\mapsto \hat v_i$, $w_i\mapsto \hat w_i$ is distance-preserving.
By construction, $\iota$ is distance-nonexpanding; therefore we only need to show that $\iota$ is distance-noncontracting.

This part is divided into two steps.

\parit{Step 1.}
In this step, we describe three groups of conditions on these 12 angles;
we show that together they guarantee that $L$ has nonnegative curvature in the sense of Alexandrov, and $\iota$ is distance-noncontracting.

First, we need to assume that the 12 angles of the triangles are at least as large as the corresponding model angles;
that is,
\[
\measuredangle \hinge {\hat x}{\hat v_i}{\hat w_j}\ge \angk{x}{v_i}{w_j}, 
\quad
\measuredangle \hinge {\hat v_i}{\hat x}{\hat w_j}\ge \angk{v_i}{x}{w_j},
\quad
\measuredangle \hinge {\hat w_j}{\hat x}{\hat v_i}\ge \angk{w_j}{x}{v_i},\leqno({*})
\]
for all $i$ and $j$.

\begin{wrapfigure}{r}{30mm}
\vskip-6mm
\centering
\includegraphics{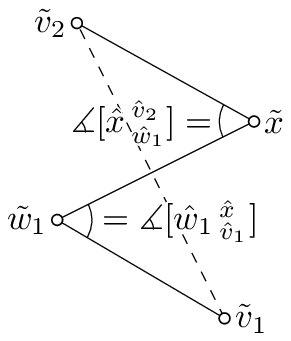}
\end{wrapfigure}

Further, choose a three-edge path in the triangulation connecting $v_1$ to $v_2$ (or $w_1$ to $w_2$), say $v_1w_1xv_2$.
Consider the plane polygonal line $\tilde v_1\tilde w_1\tilde x\tilde v_2$ with the same angles and sides as in $L$ such that $\tilde v_1$ and $\tilde v_2$ lie on the opposite sides from the line $\tilde w_1\tilde x$.
Set 
\[\tilde Z(v_1w_1xv_2)\df |\tilde v_1-\tilde v_2|.\]
The next group of conditions has eight comparisons:
\[
\begin{aligned}
|v_1- v_2|&\le \tilde Z(v_1w_ixv_2),
&
|v_1- v_2|&\le \tilde Z(v_1xw_iv_2),
\\
|w_1- w_2|&\le \tilde Z(w_1v_ixw_2),
&
|w_1- w_2|&\le \tilde Z(w_1xv_iw_2)
\end{aligned}
\leqno({*}{*})
\]
for any $i$.
Finally, we need a group of eight identities:
\[
\begin{aligned}
\measuredangle \hinge {\hat v_i}{\hat x}{\hat w_1}
+\measuredangle \hinge {\hat v_i}{\hat x}{\hat w_2}
&=\angk{v_i}{w_1}{w_2},
&
\measuredangle \hinge {\hat w_i}{\hat x}{\hat v_1}
+\measuredangle \hinge {\hat w_i}{\hat x}{\hat v_2}
&=\angk{w_i}{v_1}{v_2},
\\
\measuredangle \hinge {\hat x}{\hat v_i}{\hat w_1}
+\measuredangle \hinge {\hat x}{\hat v_i}{\hat w_2}
&=\pi,
&
\measuredangle \hinge {\hat x}{\hat w_i}{\hat v_1}
+\measuredangle \hinge {\hat x}{\hat w_i}{\hat v_2}
&=\pi
\end{aligned}
\leqno(\asterism)
\]
for any $i$.

Now, let us show that these conditions imply that $\iota$ distance-noncontracting.
Suppose $\gamma$ is a curve from $\hat x$ to $\hat v_i$ that lies completely in one of the triangles adjacent to the edge $\hat x \hat v_i$.
Note that the inequalities in $({*})$ imply that 
\[\length \gamma\ge | x -\hat v_i|_F.\]
The same holds for any pair $(\hat x,\hat w_i)$ and $(\hat v_i,\hat w_j)$.
It implies that minimizing geodesic from $x$ to any point on four edges $\hat x \hat v_i$ or $\hat x \hat w_i$ runs in the corresponding edge; in particular, we have 
\[|x- v_i|_F=|\hat x- \hat v_i|_L\qquad\text{and}\qquad |x- w_i|_F=|\hat x- \hat w_i|_L\]
for each $i$.
We also get that each of the four edges $\hat x \hat v_i$ or $\hat x \hat w_i$ is a convex set in $L$;
in particular, each of these edges can be crossed at most once by a shortest path in $L$. 

Suppose that there is a curve $\gamma$ from $\hat v_1$ to $\hat v_2$ that is shorter than $|v_1-v_2|_F$.
Since two edges $\hat v_1\hat x$ and $\hat x\hat v_2$ have total length $|v_1-v_2|_F$,
we can assume that $\gamma$ runs in a pair of two adjacent triangles, say $[\hat v_1\hat x\hat w_1]$ and $[\hat v_2\hat x\hat w_1]$.
From above, $\gamma$ crosses the edge $\hat x\hat w_1$ once.
Denote by $\hat z_1$ and $\hat z_2$ the singular points in the triangles $[\hat v_1\hat x\hat w_1]$ and $[\hat v_2\hat x\hat w_1]$.
We have the following 4 options: 

\begin{longtable}{|c|l|}
 \hline
\begin{minipage}{40mm}
\vskip3mm
\centering
\includegraphics{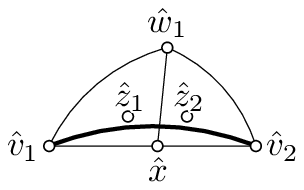}\ 
\\ \ 
\end{minipage}
&
\begin{minipage}{70mm}

\ 

If points $z_1$ and $z_2$ lie on the left from $\gamma$,
then we arrive at a contradiction with
\[\measuredangle \hinge {\hat x}{\hat w_1}{\hat v_1}
+\measuredangle \hinge {\hat x}{\hat w_1}{\hat v_2}
=\pi\]
in 
$(\asterism)$.

\ 

\end{minipage}
\\ 
\hline

\begin{minipage}{40mm}
\vskip3mm
\centering
\includegraphics{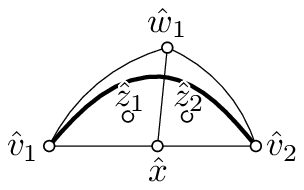}\ 
\\ \ 
\end{minipage}
&
\begin{minipage}{70mm}

\ 

If points $z_1$ and $z_2$ lie on the right from $\gamma$,
then we arrive at a contradiction with
\[
\measuredangle \hinge {\hat w_1}{\hat x}{\hat v_1}
+\measuredangle \hinge {\hat w_1}{\hat x}{\hat v_2}
=\angk{w_1}{v_1}{v_2},
\]
in 
$(\asterism)$.

\ 

\end{minipage}
\\ 
\hline

\begin{minipage}{40mm}
\vskip3mm
\centering
\includegraphics{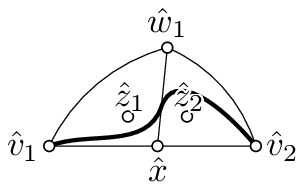}\ 
\\ \ 
\end{minipage}
&
\begin{minipage}{70mm}

\ 

If $z_1$ lies on the left side from $\gamma$, 
and $z_2$ lies on its right side, then we arrive at a contradiction with 
\[|v_1- v_2|\le \tilde Z(v_1xw_1v_2)\]
in $({*}{*})$.

\ 

\end{minipage}
\\ 
\hline

\begin{minipage}{40mm}
\vskip3mm
\centering
\includegraphics{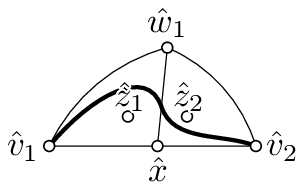}\ 
\\ \ 
\end{minipage}
&
\begin{minipage}{70mm}

\ 

If $z_1$ lies on the right side from $\gamma$, 
and $z_2$ lies on its left side, then we arrive at a contradiction with
\[|v_1- v_2|\le \tilde Z(v_1w_1xv_2)\]
in $({*}{*})$.

\ 

\end{minipage}
\\ 
\hline
\end{longtable}

It shows that $\iota$ does not decrease the distance between $v_1$ and $v_2$;
the same argument works for $w_1$ and $w_2$.
In addition, we get that two edges $\hat v_1\hat x$ and $\hat x\hat v_2$ form a shortest path in $L$; the same holds for $\hat w_1\hat x$ and $\hat x\hat w_2$.

Finally, suppose $\gamma$ is a curve from $\hat v_1$ to $\hat w_1$ that is shorter than $|v_1-w_1|_F$.
From above it does not lie in the triangle $[\hat x\hat v_1\hat w_1]$.
Recall that $\gamma$ crosses each of the four edges $\hat x \hat v_i$ or $\hat x \hat w_i$ at most once.
Therefore, $\gamma$ has to cross edge $\hat x \hat v_2$.
Since  $\hat v_1\hat x$ and $\hat x\hat v_2$ form a shortest path, we can assume that $\gamma$ visits $\hat x$ and so
\[\length \gamma\ge |v_1-x|_F+|x-w_1|_F\ge |v_1-w_1|_F\]
--- a contradiction.
The same way we show that $\iota$ does not increase the distances for each pair $(v_i,w_j)$.

It remains to show that $L$ is  Alexandrov space with nonnegative curvature.
By $(\asterism)$ the total angle around $\hat x$ in $L$ is $2\cdot\pi$.
Further, $(\asterism)$ implies that
\[
\measuredangle \hinge {\hat v_i}{\hat x}{\hat w_1}
+\measuredangle \hinge {\hat v_i}{\hat x}{\hat w_2}
\le \pi
\qquad\text{and}\qquad
\measuredangle \hinge {\hat w_i}{\hat x}{\hat v_1}
+\measuredangle \hinge {\hat w_i}{\hat x}{\hat v_2}
\le \pi.
\]
for any $i$;
that is, $L$ has convex boundary.
In particular, $L$ has locally nonnegative curvature.
It remains to apply the globalization theorem \cite[8.32]{alexander2019alexandrov}. 
(Instead, one may also apply the characterization of nonnegatively curved polyhedral spaces \cite[12.5]{alexander2019alexandrov}.)

\parit{Step 2.}
In this step we show that the 12 angles can be chosen so that they meet all the conditions $({*})$, $({*}{*})$, and $(\asterism)$.
This part is done by means of elementary geometry.

By \ref{clm:(4+1)=LSS(4)}, we can apply (4+1) comparison for the array $x$, $v_1$, $v_2$, $w_1$, $w_2$. This way we get points $\tilde x$, $\tilde v_1$, $\tilde v_2$, $\tilde w_1$, $\tilde w_2$ such that 
\begin{align*}
|\tilde x-\tilde v_i|_{\HH}&=|x-v_i|_{A},
&
|\tilde x-\tilde w_i|_{\HH}&=|x-w_i|_{A},
&
|\tilde v_i-\tilde w_j|_{\HH}&\ge |v_i-w_j|_{A},
\\
|\tilde v_1-\tilde v_2|_{\HH}&\ge |v_1-v_2|_{A},
&
|\tilde w_1-\tilde w_2|_{\HH}&\ge |w_1-w_2|_{A}.
\end{align*}
Since $v_1xv_2$ and $w_1xw_2$ are tense,
the triangle inequality implies equality in the last two inequalities;
that is, each triple of points $(\tilde v_1,\tilde x,\tilde v_2)$ and $(\tilde w_1,\tilde x,\tilde w_2)$ lies on one line.
In particular, the whole configuration lies in $\RR^2$.

Set 
\[\measuredangle \hinge {\hat x}{\hat v_i}{\hat w_j}
=
\measuredangle \hinge {\tilde x}{\tilde v_i}{\tilde w_j}\]
for all $i$ and $j$.
Since $|\tilde v_i-\tilde w_j|_{\HH}\ge |v_i-w_j|_{A}$, this choice meets four conditions in $({*})$ and the second half of the identities in $(\asterism)$.

We still need to choose the remaining 8 angles $\measuredangle \hinge {\hat v_i}{\hat w_j}{\hat x}$ and $\measuredangle \hinge {\hat w_j}{\hat v_i}{\hat x}$ for all $i$ and $j$.
To do this, we extend the configuration $\tilde x,\tilde v_1,\tilde v_2,\tilde w_1,\tilde w_2$ by 8 more points 
$\tilde v_{ij}$, $\tilde w_{ij}$ so that we can set 
\[
\measuredangle \hinge {\hat v_i}{\hat w_j}{\hat x}=\measuredangle \hinge {\tilde v_i}{\tilde w_{ij}}{\tilde x},
\qquad
\measuredangle \hinge {\hat w_j}{\hat v_i}{\hat x}=\measuredangle \hinge {\tilde w_j}{\tilde v_{ji}}{\tilde x}.
\]
We assume that $|\tilde v_i-\tilde w_{ij}|_{\RR^2}=|v_i-w_j|$ and $|\tilde w_i-\tilde v_{ij}|_{\RR^2}=|w_i-v_j|$ for all $i$ and~$j$.
The conditions $({*})$, and $({*}{*})$ will follow if we could choose the points so that
\begin{align*}
|\tilde x-\tilde v_{ij}|&\ge | x-v_j|,
&
|\tilde x-\tilde w_{ij}|&\ge | x-w_j|,
\\
|\tilde v_{j'}-\tilde v_{ij}|&\ge | v_{j'}-v_j|,
&
|\tilde w_{j'}-\tilde w_{ij}|&\ge |w_{j'}-w_j|,
\end{align*}
here we assume that $j'\ne j$, so $2'=1$ and $1'=2$.

\begin{wrapfigure}{r}{33mm}
\vskip-8mm
\centering
\includegraphics{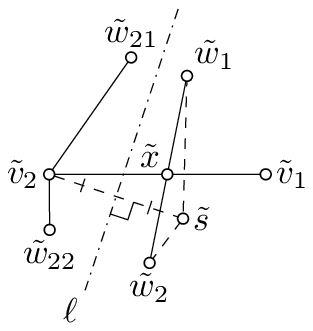}
\end{wrapfigure}

The needed points $\tilde w_{21}$ and $\tilde w_{22}$ can be chosen to be reflections of $\tilde w_{1}$ and $\tilde w_{2}$ respectively across a line $\ell$ that we are about to describe.
Suppose $[\tilde s\tilde w_1 \tilde w_2]$ is a model triangle for $[v_2w_1w_2]$ such that $\tilde s$ lies on the opposite side from $\tilde v_2$ with respect to the line $\tilde w_1\tilde w_2$.
Then $\ell$ is the perpendicular bisector of $[\tilde v_2, \tilde s]$.
Since $|\tilde w_i-\tilde v_2|\ge | w_i- v_2|=|\tilde w_i-\tilde s|$ the points $\tilde w_1$ and $\tilde w_2$ lie on the opposite side from $\tilde v_2$ with respect to $\ell$.
Whence the conditions on $\tilde w_{21}$ and $\tilde w_{22}$ follow.
By construction, we get one of the identities in $(\asterism)$ with base point $v_2$.

Similarly, we construct the remaining 6 points.

\parit{Final step.}
It remains to modify $L$ into a plane with a smooth Riemannian metric.
First, note that $L$ is a convex subset of a flat plane with at most 4 conic points.
Further, the geodesics between the 5-point subset in $L$ do not visit these conic points.
Therefore a slight smoothing around singularities does not create a problem.
\qeds

\parit{Proof of \ref{lem:key}.}
Observe that any pair of points of $F$ must lie in a tense set.
If not, then all $\LSS(4)$ inequalities will remain to hold after a slight change of the distance between the pair.
The latter contradicts that $F$ is extreme.

\begin{wrapfigure}[5]{l}{60mm}
\vskip-2mm
\centering
\includegraphics{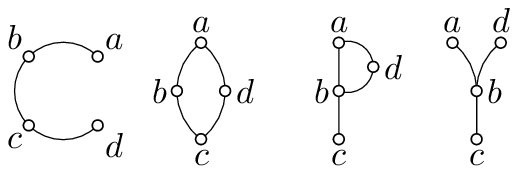}
\label{page:COPY}
\end{wrapfigure}

Suppose that two tense triples share two points.
All possible 4 configurations are shown on the diagram; they will be referred to as $C$, $O$, $P$, and $Y$ respectively.
Observe that in the configurations $P$ and $Y$, the set $\{a,b,c,d\}$ must be tense with center at~$b$.
Indeed, in the $P$-case the 4-point space is isometric to a 4-point subset on a line with order $a,d,b,c$.
In the $Y$-case, $\angk{b}{a}{c}=\angk{b}{d}{c}=\pi$, and the comparison implies that $\angk{b}{d}{a}=0$.
Without loss of generality we may assume $|a-b|\ge |d-b|$;
so, the 4-point space is also isometric to a 4-point subset on a line with order $a,d,b,c$.
That is, if $P$ or $Y$ appear in $F$, then $F$ has a 4-point tense set.
The latter contradicts the assumptions; so $P$ and $Y$ cannot appear in our configuration.

Let us show that $F$ contains at least 5 tense triples;
assume $F$ has at most $4$ of them.
By \ref{clm:10-2k}, the space of quadratic forms on $\RR^4$ contains a subspace $S$ of dimension at least 2 such that for any form $U\in S$ for all $t$ sufficiently close to zero, the forms $W+t\cdot U$ satisfy all $\LSS(4)$.
Therefore $F$ is not extremal --- a contradiction.

The remaining part of the proof is a brute-force search of all possible configurations that satisfy the conditions above.
This search is sketched on the following diagram which needs some explanation.
\begin{figure}[ht!]
\centering
\includegraphics{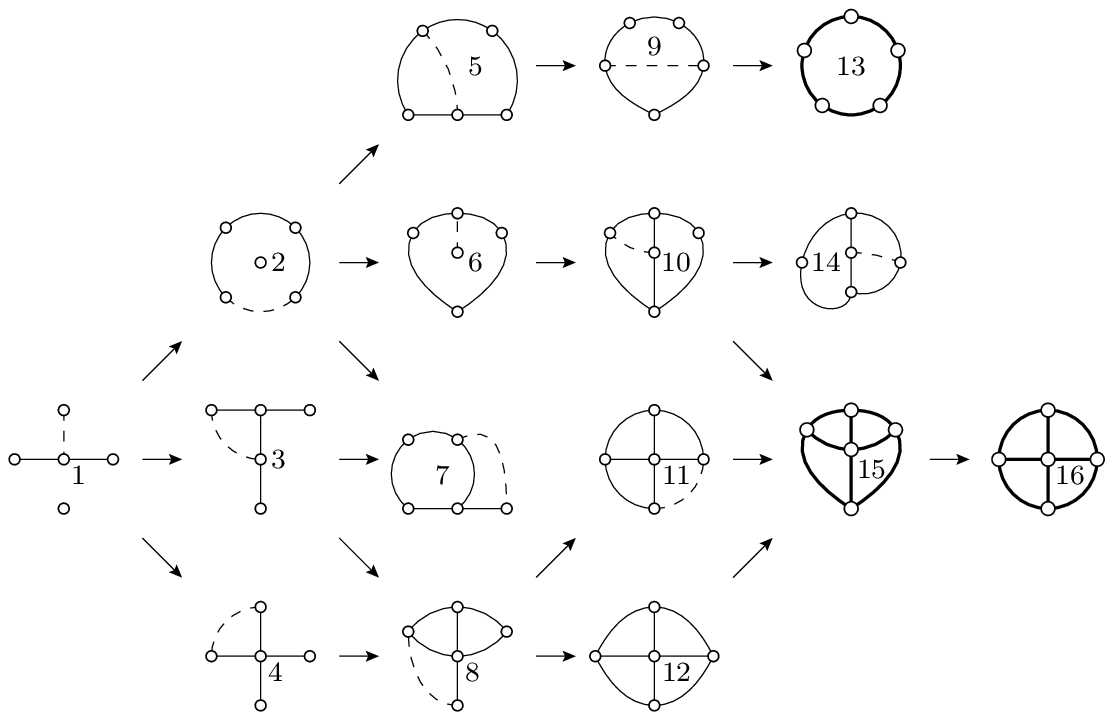}
\end{figure}\label{pic-210}
We start with a configuration with one triple marked by a solid line.
Choose a pair that is not in any triple of the configuration;
connect it by a dashed line and search for an extra triple with this pair inside.
Each time we need to check up to 9 triples that contain the pair --- 3 choices for extra points and 3 choices for the center in the obtained triple.
Some of them make a $P$ or $Y$ configuration with an existing triple, so they cannot be added.
If some of them can be added, then we draw a new diagram connected by an arrow and continue.
In many cases, symmetry reduces the number of cases.

If there are no free pairs (these are \ref{a(bcd)e+bcd}, \ref{abcdead}, \ref{abcda+aec+bed}, and \ref{abcdab+aec+bed}),
then we need to check all triples,
but due to symmetry, the number of triples can be reduced.

Once we did the classification, we need to find all configurations with at least~5 triples (these start with column 5)
such that each pair belongs to one of the tense triples (those that have no dashed line).
So we are left with three cases \ref{abcdead}, \ref{abcda+aec+bed}, and \ref{abcdab+aec+bed} marked in bold;
it proves the lemma.

The following table describes procedures at each node on the diagram.
It uses the following notations.
If a candidate triple, say $abd$ violates $Y$ rule with an existing triple, say $abc$, then we write \xcancel{$abd$}$Yabc$.
Similarly, if a candidate triple, say $adb$ violates $P$ rule with an existing triple, say $abc$, then we write \xcancel{$adb$}$Pabc$.
Further, assume a candidate, say $dbe$, does not violate the rules and so it can be added.
Suppose that after adding this triple we get a new configuration, say~\ref{abc+dbe};
in this case, we write $dbe{\to}$\ref{abc+dbe}.
Note that the new configuration is relabeled arbitrarily.

In cases \ref{a(bcd)e+bcd}, \ref{abcdead}, \ref{abcda+aec+bed}, and \ref{abcdab+aec+bed} we check all triples up to symmetry.
The used symmetries are marked in the third column.

\newcounter{foo}
\setcounter{foo}{0}
\newcommand{\myitem}{\refstepcounter{foo}\thefoo}

\begin{longtable}{|c|c|c|l|}
 \hline
\myitem\label{abc}
&
\begin{minipage}{20mm}
\vskip3mm
\centering
\includegraphics{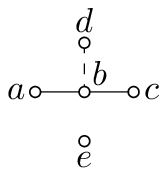}\ 
\\ \ 
\end{minipage}
&
$a\leftrightarrow c$
& 
\begin{tabular}{ll}
\xcancel{$dba$}$Yabc$;&
$dbe{\to}$\ref{abc+dbe};\\
$dab{\to}$\ref{abcd};&
$deb{\to}$\ref{abc+dae};\\
\xcancel{$adb$}$Pabc$;&
$edb{\to}$\ref{abc+dae}.\\
\end{tabular}
\\ 
\hline

\myitem\label{abcd}
&
\begin{minipage}{20mm}
\vskip3mm
\centering
\includegraphics{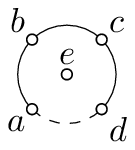}\ 
\\ \ 
\end{minipage}
&
$b\leftrightarrow c$
& 
\begin{tabular}{ll}
\xcancel{$adb$}$Pabc$;&
$ade{\to}$\ref{abcd+ade};\\
\xcancel{$abd$}$Yabc$;&
$aed{\to}$\ref{abcd+aed};\\
$bad{\to}$\ref{abcda};&
$ead{\to}$\ref{abcd+ade}.\\
\end{tabular}
\\ 
\hline

\myitem\label{abc+dae}
&
\begin{minipage}{20mm}
\vskip3mm
\centering
\includegraphics{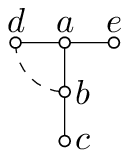}\ 
\\ \ 
\end{minipage}
&
& 
\begin{tabular}{lll}
\xcancel{$dba$}$Yabc$;&
\xcancel{$dbc$}$Yabc$;&
$dbe{\to}$\ref{abc+d(ab)e};\\
\xcancel{$dab$}$Ydae$;&
$dcb{\to}$\ref{abcd+ade};&
\xcancel{$deb$}$Pdae$;\\
\xcancel{$adb$}$Pabc$;&
\xcancel{$cdb$}$Pabc$;&
\xcancel{$edb$}$Pdae$.\\
\end{tabular}
\\ 
\hline

\myitem\label{abc+dbe}
&
\begin{minipage}{20mm}
\vskip3mm
\centering
\includegraphics{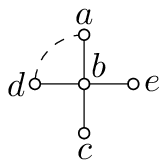}\ 
\\ \ 
\end{minipage}
&
$c\leftrightarrow e$
& 
\begin{tabular}{ll}
\xcancel{$adb$}$Pabc$;&
$adc{\to}$\ref{abc+d(ab)e};\\
\xcancel{$abd$}$Yabc$;&
\xcancel{$acd$}$Pabc$;\\
\xcancel{$bad$}$Pdbe$;&
\xcancel{$cad$}$Pabc$.\\
\end{tabular}
\\ 
\hline

\myitem\label{abcd+aed}
&
\begin{minipage}{20mm}
\vskip3mm
\centering
\includegraphics{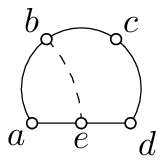}\ 
\\ \ 
\end{minipage}
&
& 
\begin{tabular}{lll}
\xcancel{$bea$}$Yaed$;&
\xcancel{$bec$}$Pabc$;&
\xcancel{$bed$}$Yaed$;\\
$bae{\to}$\ref{abcdea};&
\xcancel{$bce$}$Ybcd$;&
\xcancel{$bde$}$Pbcd$;\\
\xcancel{$abe$}$Yabc$;&
\xcancel{$cbe$}$Yabc$;&
\xcancel{$dbe$}$Pbcd$.\\
\end{tabular}
\\ 
\hline

\myitem\label{abcda}
&
\begin{minipage}{20mm}
\vskip3mm
\centering
\includegraphics{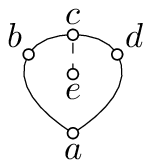}\ 
\\ \ 
\end{minipage}
&
$b\leftrightarrow d$
& 
\begin{tabular}{ll}
$cea{\to}$\ref{abcda+aec};&
\xcancel{$ceb$}$Pabc$;\\
\xcancel{$cae$}$Pabc$;&
\xcancel{$cbe$}$Yabc$;\\
\xcancel{$ace$}$Pabc$;&
\xcancel{$bce$}$Ybcd$.\\
\end{tabular}
\\ 
\hline

\myitem\label{abcd+ade}
&
\begin{minipage}{20mm}
\vskip3mm
\centering
\includegraphics{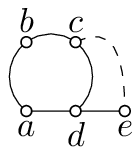}\ 
\\ \ 
\end{minipage}
&
& 
\begin{tabular}{lll}
\xcancel{$cea$}$Pade$;&
\xcancel{$ceb$}$Pabc$;&
\xcancel{$ced$}$Pbcd$;\\
\xcancel{$cae$}$Pade$;&
\xcancel{$cbe$}$Yabc$;&
\xcancel{$cde$}$Yade$;\\
\xcancel{$ace$}$Pabc$;&
\xcancel{$bce$}$Ybcd$;&
\xcancel{$dce$}$Ybcd$.\\
\end{tabular}
\\ 
\hline

\myitem\label{abc+d(ab)e}
&
\begin{minipage}{20mm}
\vskip3mm
\centering
\includegraphics{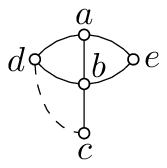}\ 
\\ \ 
\end{minipage}
&
& 
\begin{tabular}{lll}
$cda{\to}$\ref{abcd+aec+bed};&
\xcancel{$cdb$}$Pabc$;&
\xcancel{$cde$}$Pdae$;\\
\xcancel{$cad$}$Pabc$;&
\xcancel{$cbd$}$Yabc$;&
\xcancel{$ced$}$Pdae$;\\
\xcancel{$acd$}$Pabc$;&
\xcancel{$bcd$}$Pdbe$;&
$ecd{\to}$\ref{a(bcd)e+bcd};\\
\end{tabular}
\\ 
\hline

\myitem\label{abcdea}
&
\begin{minipage}{20mm}
\vskip3mm
\centering
\includegraphics{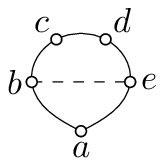}\ 
\\ \ 
\end{minipage}
&
$c\leftrightarrow d$
& 
\begin{tabular}{ll}
\xcancel{$bea$}$Ydea$;&
\xcancel{$bec$}$Pabc$;\\
$bae{\to}$\ref{abcdead};&
\xcancel{$bce$}$Ybcd$;\\
\xcancel{$abe$}$Yabc$;&
\xcancel{$cbe$}$Yabc$.\\
\end{tabular}
\\ 
\hline

\myitem\label{abcda+aec}
&
\begin{minipage}{20mm}
\vskip3mm
\centering
\includegraphics{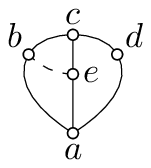}\ 
\\ \ 
\end{minipage}
&
& 
\begin{tabular}{lll}
\xcancel{$bea$}$Yaec$;&
\xcancel{$bec$}$Yaec$;&
$bed{\to}$\ref{abcda+aec+bed};\\
$bae{\to}$\ref{abcdaec};&
\xcancel{$bce$}$Ybcd$;&
\xcancel{$bde$}$Pbcd$;\\
\xcancel{$abe$}$Yabc$;&
\xcancel{$cbe$}$Yabc$;&
\xcancel{$dbe$}$Pbcd$.\\
\end{tabular}
\\ 
\hline

\myitem\label{abcd+aec+bed}
&
\begin{minipage}{20mm}
\vskip3mm
\centering
\includegraphics{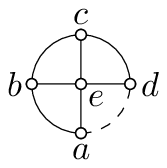}\ 
\\ \ 
\end{minipage}
&
$b\leftrightarrow c$
& 
\begin{tabular}{ll}
\xcancel{$adb$}$Pbed$;&
\xcancel{$ade$}$Paec$;\\
\xcancel{$abd$}$Yabc$;&
\xcancel{$aed$}$Ybed$;\\
$bad{\to}$\ref{abcda+aec+bed};&
\xcancel{$ead$}$Pbed$.\\
\end{tabular}
\\ 
\hline

\myitem\label{a(bcd)e+bcd}
&
\begin{minipage}{20mm}
\vskip3mm
\centering
\includegraphics{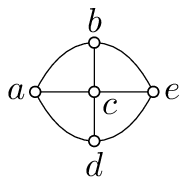}\ 
\\ \ 
\end{minipage}
&
\begin{tabular}{l}
$a\leftrightarrow e$\\
$b\leftrightarrow d$
\end{tabular}
& 
\begin{tabular}{ll}
\xcancel{$abc$}$Pace$;&
\xcancel{$abd$}$Yabe$;
\\
\xcancel{$bca$}$Yace$;&
\xcancel{$bda$}$Yade$;
\\
\xcancel{$cab$}$Pbcd$;&
$dab{\to}$\ref{abcda+aec+bed}.
\\
\end{tabular}
\\ 
\hline

\myitem\label{abcdead}
&
\begin{minipage}{20mm}
\vskip3mm
\centering
\includegraphics{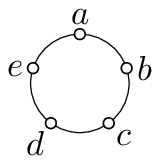}\ 
\\ \ 
\end{minipage}
&
$D_5$
& 
\begin{tabular}{l}
\xcancel{$dac$}$Pbcd$;\\
\xcancel{$acd$}$Ybcd$.\\
\end{tabular}
\\ 
\hline

\myitem\label{abcdaec}
&
\begin{minipage}{20mm}
\vskip3mm
\centering
\includegraphics{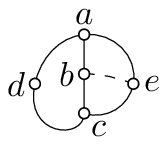}\ 
\\ \ 
\end{minipage}
&
& 
\begin{tabular}{lll}
\xcancel{$bea$}$Yaec$;&
\xcancel{$bec$}$Yaec$;&
\xcancel{$bed$}$Pdae$;\\
\xcancel{$bae$}$Ydae$;&
\xcancel{$bce$}$Ybcd$;&
\xcancel{$bde$}$Pdae$;\\
\xcancel{$abe$}$Yabc$;&
\xcancel{$cbe$}$Yabc$;&
\xcancel{$dbe$}$Pbcd$;\\
\end{tabular}
\\ 
\hline

\myitem\label{abcda+aec+bed}
&
\begin{minipage}{20mm}
\vskip3mm
\centering
\includegraphics{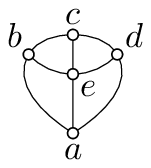}\ 
\\ \ 
\end{minipage}
&
$b\leftrightarrow d$
& 
\begin{tabular}{lll}
$bad{\to}$\ref{abcdab+aec+bed};&
\xcancel{$bec$}$Ybed$;&
\xcancel{$bea$}$Ybed$;
\\
\xcancel{$adb$}$Pbed$;&
\xcancel{$ecb$}$Ybcd$;&
\xcancel{$eab$}$Pbed$;
\\
\xcancel{$dba$}$Pbed$;&
\xcancel{$cbe$}$Yabc$;&
\xcancel{$abe$}$Yabc$;
\\
\end{tabular}
\\ 
\hline

\myitem\label{abcdab+aec+bed}
&
\begin{minipage}{20mm}
\vskip3mm
\centering
\includegraphics{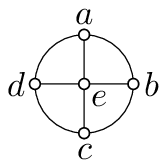}\ 
\\ \ 
\end{minipage}
&
$D_4$
& 
\begin{tabular}{l}
\xcancel{$abe$}$Yabc$;\\
\xcancel{$bea$}$Pabc$.\\
\end{tabular}
\\ 
\hline
\end{longtable}
\qedsf

\section{Final remarks}

Our theorem provides an affirmative answer to Question~6.2 in \cite{lebedeva-petrunin}.
For 6-point metric spaces, a direct analog of the theorem does not hold,
but Question~6.3 in \cite{lebedeva-petrunin} contains the corresponding conjecture.

An analogous problem for 5-point sets in nonpositively curved spaces was solved by Tetsu Toyoda \cite{toyoda}; another solution is given in \cite{lebedeva-petrunin}.
The 6-point case is open;
see  \cite[Question 6.1]{lebedeva-petrunin} and a partial answer in \cite{lebedeva-petrunin-octahedron}.

{

\begin{wrapfigure}{r}{30mm}
\vskip-3mm
\centering
\includegraphics{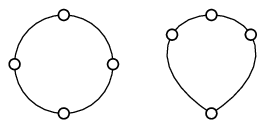}
\end{wrapfigure}

The 4-point case is much easier;
the classification gives only two cases on the diagram.
Both admit an embedding into a circle.
It can be used to prove the following statement.

}

\begin{thm}{Theorem}
Any 4-point space satisfying the nonnegative-curvature comparison
admits an embedding into a product of a circle and Euclidean space.
\end{thm}

The following statement can be proved in a similar manner.

\begin{thm}{Theorem} Any 4-point space satisfying the nonpositive-curvature comparison admits an embedding into a product of a tripod and Euclidean space.
\end{thm}

These two results are analogous to Wald's theorem mentioned in the introduction;
they were obtained in the note by Vladimir Zolotov and the first author~\cite{lebedeva-zolotov}.

\medskip

Most of our arguments can be applied to arbitrary curvature bound;
Section~\ref{sec:ext} is the only place where we essentially use that the bound is zero.

It would be interesting to classify 5-point subsets in other classes of spaces;
for example, in products of circles, or complete flat manifolds.
Note that the second and third types of spaces in the classification lemma (\ref{lem:key}) do not admit an embedding into a product of circles; so the answer must be different.
We are not aware of 5-point spaces that admit an isometric embedding into a complete nonnegatively curved Riemannian manifold, but not in a complete flat manifold.
Analogous questions can be asked about products of trees \cite{lebedeva-petrunin-octahedron} and Euclidean buildings;
these are especially nice classes of spaces with nonpositive curvature.

Our argument can be applied to attack the 6-point case \cite[Question~6.3]{lebedeva-petrunin}, except we could no longer use $\LSS(5)$.
The case of extremal metrics only with 3-point tense sets seems to be easier. 
The 5-point tense set can be solved following our argument in \ref{prop:4-tense}.
At the moment we do not see a way to do the 6- and 4-point tense sets.

It might be interesting to find conditions on finite subsets of metric spaces that are related to other curvature bounds as, for example, nonnegative curvature operator or nonnegative isotropic curvature; see \cite[1.19$_+(e)$]{gromov}.
According to \cite{lebedeva-petrunin-graph-alex}, graph comparison can be used to describe conditions that are stronger than nonnegative or nonpositive in the sense of Alexandrov, but nonnegative curvature operator has a chance to be described this way.

\switch{}{
\section*{Compliance with Ethical Standards}
The authors have declared no conflicts of interest for this article.

Data sharing not applicable to this article as no datasets were generated or analysed during the current study.
}

{\sloppy
\printbibliography[heading=bibintoc]

@book{alexander2019alexandrov,
      title={Alexandrov geometry: foundations}, 
      author={Alexander, S. and Kapovitch, V. and Petrunin, A.},
      year={2022},
      eprint={1903.08539},
      archivePrefix={arXiv},
      primaryClass={math.DG}
}

@inproceedings {AKP-Kirszbraun,
    AUTHOR = {Alexander, S. and Kapovitch, V. and Petrunin, A.},
     TITLE = {Alexandrov meets {K}irszbraun},
 BOOKTITLE = {Proceedings of the {G}\"{o}kova {G}eometry-{T}opology {C}onference
              2010},
     PAGES = {88--109},
 PUBLISHER = {Int. Press, Somerville, MA},
      YEAR = {2011},
   MRCLASS = {53C23},
  MRNUMBER = {2931882},
MRREVIEWER = {Koichi Nagano},
}

@book {gromov,
    AUTHOR = {Gromov, M.},
     TITLE = {Metric structures for {R}iemannian and non-{R}iemannian
              spaces},
    SERIES = {Modern Birkh\"{a}user Classics},
 %PUBLISHER = {Birkh\"{a}user Boston, Inc., Boston, MA},
      YEAR = {2007},
    % PAGES = {xx+585},
      ISBN = {978-0-8176-4582-3; 0-8176-4582-9},
   MRCLASS = {53C23 (53-02)},
  MRNUMBER = {2307192},
}

@article{toyoda,
    AUTHOR = {Toyoda, T.},
     TITLE = {An intrinsic characterization of five points in a
              {CAT}(0) space},
   JOURNAL = {Anal. Geom. Metr. Spaces},
  FJOURNAL = {Analysis and Geometry in Metric Spaces},
    VOLUME = {8},
      YEAR = {2020},
    NUMBER = {1},
     PAGES = {114--165},
   MRCLASS = {53C23 (51F99)},
  MRNUMBER = {4141370},
       DOI = {10.1515/agms-2020-0111},
       URL = {https://doi.org/10.1515/agms-2020-0111},
}

@misc{lebedeva-petrunin-octahedron,
  doi = {10.48550/ARXIV.2212.06445},
  url = {https://arxiv.org/abs/2212.06445},
  author = {Lebedeva, N. and Petrunin, A.},  
  title = {Trees meet octahedron comparison},  
  eprint={2212.06445},
  archivePrefix={arXiv},
  primaryClass={math.MG},
  year = {2022},
copyright = {Creative Commons Zero v1.0 Universal}
}

@misc{lebedeva-petrunin-graph-alex,
  doi = {10.48550/ARXIV.2212.08016},  
  url = {https://arxiv.org/abs/2212.08016},
  author = {Lebedeva, N. and Petrunin, A.},
  title = {Graph comparison meets Alexandrov},
  eprint={2212.08016},
  archivePrefix={arXiv},
  primaryClass={math.MG},
  year = {2022},
  copyright = {Creative Commons Zero v1.0 Universal}
}

@article {lang-schroeder,
    AUTHOR = {Lang, U. and Schroeder, V.},
     TITLE = {Kirszbraun's theorem and metric spaces of bounded curvature},
   JOURNAL = {Geom. Funct. Anal.},
  FJOURNAL = {Geometric and Functional Analysis},
    VOLUME = {7},
      YEAR = {1997},
    NUMBER = {3},
     PAGES = {535--560},
      ISSN = {1016-443X},
   MRCLASS = {53C23},
  MRNUMBER = {1466337},
MRREVIEWER = {Zhongmin Shen},
       DOI = {10.1007/s000390050018},
       URL = {https://doi.org/10.1007/s000390050018},
}

@article {lebedeva-petrunin,
    AUTHOR = {Lebedeva, N. and Petrunin, A.},
     TITLE = {5-point {$\rm CAT(0)$} spaces after {T}etsu {T}oyoda},
   JOURNAL = {Anal. Geom. Metr. Spaces},
  FJOURNAL = {Analysis and Geometry in Metric Spaces},
    VOLUME = {9},
      YEAR = {2021},
    NUMBER = {1},
     PAGES = {160--166},
   MRCLASS = {53C23 (51F99)},
  MRNUMBER = {4298102},
       DOI = {10.1515/agms-2020-0126},
       URL = {https://doi.org/10.1515/agms-2020-0126},
}

@article {lebedeva-petrunin-zolotov,
    AUTHOR = {Lebedeva, N. and Petrunin, A. and Zolotov, V.},
     TITLE = {Bipolar comparison},
   JOURNAL = {Geom. Funct. Anal.},
  FJOURNAL = {Geometric and Functional Analysis},
    VOLUME = {29},
      YEAR = {2019},
    NUMBER = {1},
     PAGES = {258--282},
      ISSN = {1016-443X},
   MRCLASS = {53C23},
  MRNUMBER = {3925110},
MRREVIEWER = {Nan Li},
       DOI = {10.1007/s00039-019-00481-9},
       URL = {https://doi.org/10.1007/s00039-019-00481-9},
}

@online{lebedeva-zolotov,
author = {Lebedeva, N. and Zolotov, V.},
title = {Curvature and 4-point subspaces},
eprint ={www.researchgate.net/publication/367511123_Curvature_and_4-point_subspaces}
}

@MISC{perelman:spaces2,
AUTHOR = {Perelman, G.},
    TITLE = {Alexandrov paces with curvature
bounded from below II.},
HOWPUBLISHED = {Preprint LOMI, 1991.},
    EPRINT = {https://anton-petrunin.github.io/papers/},
}

@article {petrunin-2017,
AUTHOR = {Petrunin, A.},
TITLE = {In search of a five-point Alexandrov type condition},
JOURNAL = {St. Petersburg Math. J.},
VOLUME = {29},
YEAR = {2018},
NUMBER = {1},
PAGES = {223--225},
}

@article {sturm,
    AUTHOR = {Sturm, K.-T.},
     TITLE = {Metric spaces of lower bounded curvature},
   JOURNAL = {Exposition. Math.},
  FJOURNAL = {Expositiones Mathematicae. International Journal},
    VOLUME = {17},
      YEAR = {1999},
    NUMBER = {1},
     PAGES = {35--47},
      ISSN = {0723-0869},
   MRCLASS = {53C23},
  MRNUMBER = {1687468},
MRREVIEWER = {Igor Belegradek},
}

@article{wald,
  title={Begr\"undung eiiner Koordinatenlosen Differentialgeometrie der Fl\"achen},
  author={Wald, A.},
  journal={Ergebnisse eines mathematischen Kolloquium},
  volume={6},
  pages={24--46},
  year={1935},
language={german},
hyphenation={german}
}
\fussy
}

\Addresses
\end{document}